\documentclass[12pt,a4paper]{amsart}

\usepackage{graphicx}
\usepackage{geometry}
\begin{document}

\noindent
{\LARGE\bf Pick's Theorem in Two-Dimensional Subspace of $\mathbb{R}^3$

}

\vspace{0.5cm} \noindent{\large Lin Si}

\noindent
College of Science, Beijing Forestry University, Beijing, 100083, P.R.China

\noindent
E-mail: silincd@163.com

\vspace{0.5cm}
\noindent
{\bf Abstract.} In this note, we given a version of Pick's theorem for the simple lattice polygon in two-dimensional subspace of $\mathbb{R}^3$.



\vspace{0.9cm}\noindent
{\large\bf 1. Introduction}

\vspace{0.4cm}\noindent

In the Euclidean plane $\mathbb{R}^2$, a lattice point is one whose coordinates are both integers. A lattice polygon is a polygon with all vertices on integer coordinates. The area $A(P)$ of a simple lattice polygon $P$ can be given by the celebrated Pick's theorem[5]

$$A(P)=I(P)+\frac{1}{2}B(P)-1,$$
where $B(P)$ is the number of lattice points on the boundary of $P$ and $I(P)$ is the number of lattice points in the interior of $P$.

Pick's formula can be used to compute the area of a lattice polygon conveniently.

For example, on Figure 1, $I(P)=60$, $B(P)= 15$. Then, the area of the polygon is
$A(P) =60 + 15 - 1 = 74$.


\begin{figure}[htp]
\includegraphics[height=6cm,width=8cm,angle=0]{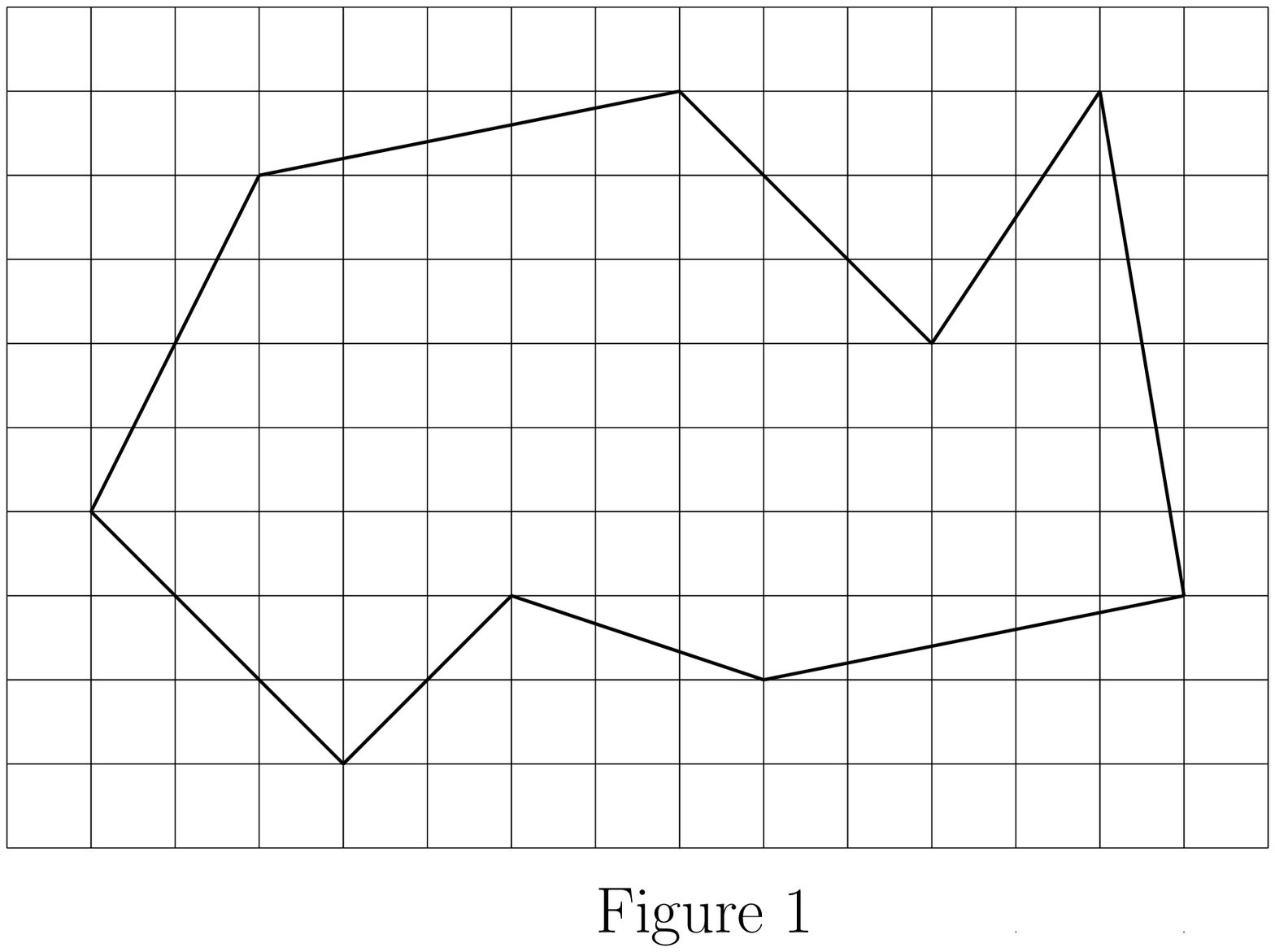}
\end{figure}

There are many papers concerning Pick's theorem and it's generalizations[1,2,3,4], which mostly be discussed in two dimension.

Unfortunately, Pick theorem is failed in three dimensions.  In 1957, John Reeve found a class of tetrahedra, named as Reeve tetrahedra later, whose
vertices are
$$(0,0,0)^T, (1,0,0)^T, (0,1,0)^T,  (1,1,r)^T,$$
where $r$ is a positive integer.

All Reeve tetrahedra contain the same number of lattice points, but their volumes are different.

In this note, we discussed  Pick's theorem in two-dimensional subspace of $\mathbb{R}^3$.  For any $(a,b,c)^T\in \mathbb{Z}^3$ with $(a,b,c)=1$, i.e. the greatest common factor of $a,b,c$ is one, denote
by $K$, $ax+by+cz=0$, the two-dimensional subspace of $\mathbb{R}^3$. Then we established the following theorem.

\medskip
\noindent
{\bf Theorem.} {\it If $P$ is simple lattice polygon in the K, then the area of $P$ is
$$k(I(P)+\frac{1}{2}B(P)-1),$$
where $B(P)$ is the number of lattice points on the boundary of $P$ in $\mathbb{Z}^3$, $I(P)$ is the number of lattice points in the interior of $P$ in $\mathbb{Z}^3$
and $k$ is the constant $(a^3+ab^2)\sqrt{a^2+b^2+c^2}$. }

\medskip
\noindent
{\bf Remark.} {\it Although the simple lattice polygon $P$ is in the two-dimensional subspace $K$, the lattice points in $P$ belong to $\mathbb{Z}^3$.}

Let $(a,b,c)^T=(1,0,0)^T$ in the Theorem, then we can get Pick's theorem in some coordinate plane of $\mathbb{R}^3$.

\medskip
\noindent
{\bf Corollary.} {\it If $P$ is simple lattice polygon in the $K$, whose normal vector is $(1,0,0)^T$, then the area of $P$ is
$$I(P)+\frac{1}{2}B(P)-1.$$}

\vspace{0.9cm}\noindent
{\large\bf 2. Proof of Main Result}

 For any $(a,b,c)^T\in \mathbb{Z}^3$ with $(a,b,c)=1$, there is two-dimensional subspace of $\mathbb{R}^3$
$$ax+by+cz=0,\eqno (1)$$
whose normal vector is just  $(a,b,c)^T$. We denote this two-dimensional subspace by $K$.

By the theory of linear equations system,  $(-b,a,0)^T$ and $(-c,0,a)^T$ are two linearly independent solutions of (1).
We denote $(-b,a,0)^T$ by  $\alpha$ and $(-c,0,a)^T$ by $\beta$. 	
Obviously, $\alpha$ and $\beta$ are also the basis of $K$.

\medskip
\noindent
{\bf Lemma 1.} {\it For any $(a,b,c)^T\in \mathbb{Z}^3$ with $(a,b,c)=1$, there exists the lattice basis with the minimal  area in the two-dimensional subspace K.}

\medskip\noindent
{\bf Proof.}  The area of parallelogram generated by  $\alpha$ and $\beta$ is

$$\frac{1}{\sqrt{a^2+b^2+c^2}}\begin{vmatrix}
a & -b &  -c\\
b & a  &   0\\
c & 0  &   a\\
\end{vmatrix}=\frac{a^3+ab^2+ac^2}{\sqrt{a^2+b^2+c^2}}=a\sqrt{a^2+b^2+c^2}.$$

Denote $(a,b,c)^T$ by $n$.
For any lattice basis in $K$, $k_1\alpha+k_2\beta$ and $l_1\alpha+l_2\beta$, where $k_i, l_i\in \mathbb{Z}(i=1,2)$ and
$\begin{vmatrix}
k_1 &  l_1 \\
k_2 &  l_2  \\
\end{vmatrix}=0$.
The area of parallelogram generated by  $k_1\alpha+k_2\beta$ and $l_1\alpha+l_2\beta$ is

$$\frac{1}{\sqrt{a^2+b^2+c^2}}\begin{vmatrix}
n & \vdots & k_1\alpha+k_2\beta & \vdots &  l_1\alpha+l_2\beta\\
\end{vmatrix}=\frac{1}{\sqrt{a^2+b^2+c^2}}\begin{vmatrix}
(n, \alpha, \beta) \times\left(\begin{matrix}
1    &    0    &    0\\
0    &   k_1 &  l_1 \\
0    &   k_2 &  l_2 \\
\end{matrix}\right)\\
\end{vmatrix},$$
where $\begin{vmatrix}
n & \vdots & k_1\alpha+k_2\beta & \vdots &  l_1\alpha+l_2\beta\\
\end{vmatrix}$ denote the determinant of $n, k_1\alpha+k_2\beta$ and $l_1\alpha+l_2\beta$.

Thus the lattice basis $k_1\alpha+k_2\beta$ and $l_1\alpha+l_2\beta$ have the minimal  area if and only if $\begin{vmatrix}
k_1 &  l_1 \\
k_2 &  l_2  \\
\end{vmatrix}=1$.

Let $k_1=1, k_2=0,  l_1=0, l_2=1$, $\alpha, \beta$ are the lattice basis with the minimal  area in the two-dimensional subspace $K$. \hfill{$\Box$}

\medskip
\noindent
{\bf Lemma 2.} {\it For any $(a,b,c)^T\in \mathbb{Z}^3$ with $(a,b,c)=1$, there exists the orthogonal lattice basis in the two-dimensional subspace K.}

\medskip\noindent
{\bf Proof.}  By Lemma 1, $\alpha, \beta$ are the lattice basis with the minimal  area in the two-dimensional subspace $K$.
By Schmidt orthogonalization, let

\begin{eqnarray*}
\hspace{3.6cm}
\gamma_1&\hspace{-0.2cm}=&\hspace{-0.2cm}\alpha=(-b,a,0)^T\\ 
\gamma_2  &\hspace{-0.2cm}=&\hspace{-0.2cm}\beta-\frac{(\beta,\gamma_1)}{(\gamma_1,\gamma_1)}\\ 
&\hspace{-0.2cm}=&\hspace{-0.2cm}(-c,0,a)^T-\frac{bc}{a^2+b^2}\gamma_1\\
&\hspace{-0.2cm}=&\hspace{-0.2cm}(-c,0,a)^T-(\frac{-b^2c}{a^2+b^2},\frac{abc}{a^2+b^2},0)^T\\
&\hspace{-0.2cm}=&\hspace{-0.2cm}(\frac{-a^2c}{a^2+b^2},\frac{-abc}{a^2+b^2},\frac{a^3+ab^2}{a^2+b^2})^T,\hspace{6.8cm}
\end{eqnarray*}
where $(\beta,\gamma_1)$ denote the usual inner product of $\beta,\gamma_1$ in $\mathbb{R}^3$.

Thus 
\begin{eqnarray*}
\eta_1=\gamma_1&\hspace{-0.2cm}=&\hspace{-0.2cm}\alpha=(-b,a,0)^T\\
\eta_2=(a^2+b^2)\gamma_2  &\hspace{-0.2cm}=&\hspace{-0.2cm}(-a^2c,-abc,a^3+ab^2)^T
\end{eqnarray*}
are the orthogonal lattice basis in the two-dimensional subspace K.\hfill{$\Box$}

\medskip
\noindent
{\bf Proof of Theorem} By Lemma 2, $\eta_1, \eta_2$ are the orthogonal lattice basis in the two-dimensional subspace $K$.

The area of parallelogram generated by  $\eta_1, \eta_2$ is

\begin{eqnarray*}
\hspace{3.6cm}
   &\hspace{-0.2cm}  &\hspace{-0.2cm}\frac{1}{\sqrt{a^2+b^2+c^2}}\begin{vmatrix}
a & -b &  -a^2c\\
b & a  &   -abc\\
c & 0  &   a^3+ab^2\\
\end{vmatrix}\\
&\hspace{-0.2cm}=&\hspace{-0.2cm}\frac{a^5+a^3b^2+ab^2c^2+a^3c^2+a^3b^2+ab^4}{\sqrt{a^2+b^2+c^2}}\\
&\hspace{-0.2cm}=&\hspace{-0.2cm}\frac{a(a^4+2a^2b^2+b^4+(a^2+b^2)c^2)}{\sqrt{a^2+b^2+c^2}}\\
&\hspace{-0.2cm}=&\hspace{-0.2cm}\frac{a(a^2+b^2)(a^2+b^2+c^2)}{\sqrt{a^2+b^2+c^2}}\\
&\hspace{-0.2cm}=&\hspace{-0.2cm}(a^3+ab^2)\sqrt{a^2+b^2+c^2},\hspace{6.8cm}
\end{eqnarray*}
which just is the constant $k$ in the Theorem.
\hfill{$\Box$}

\vspace{0.5cm}\noindent
{\bf Acknowledgements.}

This work was supported by the Beijing Higher Education Young Elite Teacher Project (Granted No.YETP0770) , the National Natural Science
Foundation of China (Grant No.11001014), and the Young Teachers Domestic Visiting Scholars Program of Beijing Forestry University.


\begin{thebibliography}{99}

\bibitem{ding}Ding R, Kolodziejczyk K, Reay J,  A new Pick-type theorem on the hexagonal lattice, {\it  Discrete Math.} {\bf68} (1988)171-177


\bibitem{funk} Funkenbusch W W, From Euler's formula to Pick's formula using an edge theorem, {\it The Amer. Math. Monthly} {\bf 81} (1974)647-648


\bibitem{gs}Grunbaum B, Shephard G C, Pick's theorem, {\it The Amer. Math. Monthly}  {\bf 100} (1993) 150-161


\bibitem{liu}Liu A C F, Lattice points and Pick's theorem, {\it  Math. Mag.}  {\bf52} (1979)232-235


\bibitem{pick}Pick G A,  Geometrisches zur Zahlentheorie, Sitzungber. Lotos, Naturwissen Zeitschrift  {\bf 19 }(1899) 311-319





\end{thebibliography}
\end{document}